\documentclass[a4paper]{amsart}
\usepackage{color}
\usepackage[latin1]{inputenc}
\usepackage[T1]{fontenc}
\usepackage{amsfonts}
\usepackage{amssymb}
\usepackage{amsmath}
\usepackage{amsthm}
\usepackage{graphicx,type1cm,eso-pic,color}
\usepackage{pstricks,pst-plot,pstricks-add}
\usepackage{float}
\newtheorem{theorem}{Theorem}[section]
\newtheorem{lemma}[theorem]{Lemma}
\newtheorem{proposition}[theorem]{Proposition}
\newtheorem{corollary}[theorem]{Corollary}
\newtheorem{definition}[theorem]{Definition}

\begin{document}
\setlength\arraycolsep{2pt}
\title{A Pathwise Fractional one Compartment Intra-Veinous Bolus Model}
\author{Nicolas MARIE}
\address{Laboratoire Modal'X Universit\'e Paris-Ouest 92000 Nanterre}
\email{nmarie@u-paris10.fr}
\date{}
\maketitle
%

% Abstract.

%
\begin{abstract}
Extending deterministic compartments pharmacokinetic models as diffusions seems not realistic on biological side because paths of these stochastic processes are not smooth enough. In order to extend one compartment intra-veinous bolus models, this paper suggests to model the concentration process $C$ by a class of stochastic differential equations driven by a fractional Brownian motion of Hurst parameter belonging to $]1/2,1[$.
\\
The first part of the paper provides probabilistic and statistical results on the concentration process $C$ : the distribution of $C$, a control of the uniform distance between $C$ and the solution of the associated ordinary differential equation, and consistent estimators of the elimination constant, of the driving signal's Hurst parameter and of the volatility constant.
\\
The second part of the paper provides applications of these theoretical results on simulated concentrations : a qualitative procedure to choose parameters on small sets of observations, and simulations of estimators of the elimination constant and of the driving signal's Hurst parameter. The relationship between the estimations quality and the size/length of the sample is discussed.
\end{abstract}
\tableofcontents
\noindent
%\textbf{MSC2010 :}
%\\
%\\
\textit{Corresponding author}: Nicolas MARIE.\\
\textit{Address}: Laboratoire Modal'X. Universit\'e Paris 10.\\
\indent 200 Avenue de la R\'epublique, 92000 Nanterre, France.\\
\textit{E-mail address}: nmarie@u-paris10.fr (N. MARIE).\\
\textit{Tel.}: (+33)7 70 01 83 51.
\\
\\
\textbf{Acknowledgements.} Many thanks to Francis Lavergne M.D. for his advices about possible clinical applications of that paper's results.
%

% Section : Introduction.

%
\section{Introduction}
\noindent
Compartments pharmacokinetic models describe the way an administered drug is transmitted among the \textit{body's compartments}. The concentration of the drug in each compartment can be modeled by ordinary differential equations (cf. Y. Jacomet \cite{JACOMET89}).
\\
In particular, in one compartment models, the concentration is classically modeled by a linear (deterministic) differential equation with negative constant coefficient, taking in account the absorption and elimination steps. Only one compartment models are studied in this paper.
\\
By D. D'Argenio and K. Park \cite{DP97}, the elimination process has both deterministic and random components. A natural way to take in account these components is to add a stochastic noise in the linear differential equation that classically models the concentration. That has been well studied in the It\^o stochastic calculus framework by many authors (cf. S. Donnet and A. Samson \cite{DS13}).
\\
However, as mentioned in M. Delattre and M. Lavielle \cite{DL11}, since the standard Brownian motion has $\alpha$-H\"older continuous paths with $\alpha\in ]0,1/2[$, the extension of the deterministic model as a diffusion is not realistic on biological side. M. Delattre and M. Lavielle \textit{force} the paths regularity of the concentration process $C$ by putting
\begin{displaymath}
C_t :=
C_0\exp\left(-\int_{0}^{t}D_sds\right)
\textrm{ ; }
t\in\mathbb R_+
\end{displaymath}
where $D$ is the diffusion that extends the deterministic model.
\\
As mentioned in N. Marie \cite{MARIE12}, another way to increase the regularity of the concentration process paths is to replace the standard Brownian motion by a fractional Brownian motion $B^H$ of Hurst parameter $H\in ]1/2,1[$ as driving signal. Since the signal is not a semi-martingale anymore, the stochastic integral is taken pathwise, in the sense of Young (cf. A. Lejay \cite{LEJAY10}). The Young integral keeps the regularity of the driving signal, therefore the concentration process has $\alpha$-H\"older continuous paths with $\alpha\in ]0,H[$.
\\
In both It\^o and pathwise stochastic calculus frameworks, an interesting volatility function is $x\in\mathbb R_+\mapsto\sigma x^{\beta}$ with $\sigma\in\mathbb R$ and $\beta\in [0,1]$. It covers classical models :
\begin{itemize}
 \item $\beta = 0$, $\sigma\not= 0$ : Langevin equation. Its solution is the so-called Ornstein-Uhlenbeck process.
 \item $\beta = 1/2$, $\sigma\not= 0$ : Cox-Ingersoll-Ross model.
 \item $\beta = 1$, $\sigma\not= 0$ : Linear stochastic differential equation.
 \item $\sigma = 0$ : Linear ordinary differential equation.
\end{itemize}
In the It\^o stochastic calculus framework, that concentration model has been studied on statistical side in K. Kalogeropoulos et al. \cite{KDP08}. In the pathwise stochastic calculus framework, it has been studied on probabilistic side in N. Marie \cite{MARIE12}.
\\
\\
This paper is devoted to the probabilistic and statistical study of the special case of the one compartment intra-veinous (i.v.) bolus model with fractional Brownian signal :
\begin{equation}\label{main_model}
C_t = C_0 -\upsilon\int_{0}^{t}C_sds +\sigma\int_{0}^{t}C_{s}^{\beta}dB_{s}^{H}
\textit{ $;$ }
t\in [0,\tau_0]
\end{equation}
where
\begin{displaymath}
\tau_0 :=
\inf\left\{t\in\mathbb R_+ : C_t = 0\right\},
\end{displaymath}
the exponent $\beta$ belongs to $[0,1[$, $\upsilon > 0$ is the rate of elimination describing the removal of the drug by all elimination processes including excretion and metabolism, and $C_0 := A_0/V$ with $A_0 > 0$ the administered dose and $V > 0$ the volume of the elimination compartment.
\\
Since its vector field is $C^{\infty}$ on bounded sets of $\mathbb R_{+}^{*}$, equation (\ref{main_model}) admits a unique continuous pathwise solution defined on $[0,\tau_0]$ and satisfying $C_. = X_{.}^{\gamma + 1}$, where $\gamma :=\beta/(1-\beta)$ and $X$ is the solution of the following fractional Langevin \mbox{equation :}
\begin{displaymath}
X_t =
C_{0}^{1-\beta} -
\upsilon(1-\beta)
\int_{0}^{t}X_sds +
\sigma(1-\beta)B_{t}^{H}
\textrm{ ; }
t\in\mathbb R_+.
\end{displaymath}
That equation is obtained by applying the \textit{rough} change of variable formula to the process $C$ and to the map $x\in\mathbb R_+\mapsto x^{1-\beta}$ on $[0,\tau_0]$. For details, the reader can refer to N. Marie \cite{MARIE12}. The fractional Langevin equation is deeply studied in P. Cheridito et al. \cite{CKM03}.
\\
\\
Since the concentration process has to be positive and stop when it hits zero, it can be defined as the solution of equation (\ref{main_model}) on $[0,\tau_0]$.
\\
For the sake of simplicity, even if the following equality only holds on $[0,\tau_0]$, throughout this paper, $C$ is defined on $\mathbb R_+$ by
\begin{displaymath}
C_t :=
\left|C_{0}^{1-\beta} +
\sigma B_{t}^{H}(\vartheta)
\right|^{\gamma + 1}
e^{-\upsilon t}
\textrm{ ; }
t\in\mathbb R_+
\end{displaymath}
with
\begin{displaymath}
\vartheta_t :=
(1-\beta)
e^{\upsilon(1-\beta)t}
\textrm{ ; }
t\in\mathbb R_+
\end{displaymath}
and the Young/Wiener integral (cf. Appendix A)
\begin{displaymath}
B_{t}^{H}(\vartheta) :=
\int_{0}^{t}\vartheta_sdB_{s}^{H}
\textrm{ ; }
t\in\mathbb R_+.
\end{displaymath}
Note that for $H = 1$ and $\beta = 0$, the fractional Brownian motion is matching with $t\in\mathbb R_+\mapsto \xi t$ such that $\xi\rightsquigarrow\mathcal N(0,1)$, and
\begin{displaymath}
C_t =
\left|\frac{\sigma}{\upsilon}\xi +
\left(C_{0} - \frac{\sigma}{\upsilon}\xi\right)e^{-\upsilon t}\right|
\textrm{ ; }
t\in\mathbb R_+.
\end{displaymath}
That \textit{limit} case illustrates that the Hurst parameter $H$ is continuously controlling the regularity of the concentration process paths, but also that $\sigma$ and $H$ provides two complementary ways to control the impact of the random component on the elimination process with respect to its deterministic component.
\\
\\
In mathematical finance, the semi-martingale property of the prices process is crucial in order to ensure the market's completeness. The It\^o stochastic calculus is then tailor-made to model prices in finance. In pharmacokinetic, the semi-martingale property of the concentration process seems not crucial on biological side.
\\
To replace the standard Brownian motion by a fractional Brownian motion in the pathwise stochastic calculus framework implies that the concentration process doesn't satisfy the Markov property anymore. In general, it makes the estimation of parameters $\upsilon$, $\sigma$ and $H$ difficult, but the relationship between $C$ and $X$ mentioned above allows to bypass these difficulties by using results coming from Y. Hu and D. Nualart \cite{HN10}, J. Istas and G. Lang \cite{IL97} and, A. Brouste and S. Iacus \cite{BI13}.
\\
\\
The second section is devoted to probabilistic and statistical properties of processes $X$ and $C$. The first part concerns the distribution of the concentration process $C$ and a control, in probability, of the uniform distance between the fractional Ornstein-Uhlenbeck process $X$ and the solution of the associated ordinary differential equation. The second part provides a strongly consistent estimator of the elimination constant $\upsilon$, and an extension of existing ergodic theorems for the fractional Ornstein-Uhlenbeck process $X$ is established. The third part provides a strongly consistent estimator of $(H,\sigma)$. A weakly consistent estimator of $\upsilon$ is deduced for unknown values of $H$ and $\sigma$.
\\
The third section is devoted to the application of the second subsection's theoretical results on simulated concentrations. For small sets of observations, the first part provides a qualitative procedure for choosing parameters $H$, $\sigma$ and $\beta$. The cornerstone of the procedure is the control of the uniform distance between $X$ and the solution of the associated ordinary differential equation mentioned above. The second part illustrates the convergence of estimators provided at Section 2. The relationship between the estimations quality and the size/length of the sample is discussed.
\\
Appendices A and B provide respectively useful definitions and results on fractional Brownian motions, and proofs of results stated at Section 2.
%

% Section : Probabilistic and statistical properties of the concentration process.

%
\section{Probabilistic and statistical properties of the concentration process}
\noindent
The first subsection concerns the distribution of the concentration process $C$ by using that $C_t = |X_t|^{\gamma + 1}$ ; $t\in\mathbb R_+$. Lemma \ref{auxiliary_distributions} provides the covariance function of the fractional Ornstein-Uhlenbeck process $X$. Proposition \ref{variance_upper_bound} allows to control, in probability, the uniform distance between the process $X$ and the solution of the associated ordinary differential equation. Refer to Appendix B for proofs of these results.
\\
Essentially inspired by Y. Hu and D. Nualart \cite{HN10}, the second subsection provides an ergodic theorem for the concentration process and its application to the estimation of the parameter $\upsilon$. An extension of existing ergodic theorems for the fractional Ornstein-Uhlenbeck process $X$ is established.
\\
Essentially inspired by J. Istas and G. Lang \cite{IL97} and, A. Brouste and S. Iacus \cite{BI13}, the third subsection provides a strongly consistent estimator of the parameter $(H,\sigma)$. A weakly consistent estimator of $\upsilon$ is deduced for unknown values of $H$ and $\sigma$.
%

% Subsection : Distribution of the concentration process and related topics.

%
\subsection{Distribution of the concentration process and related topics}
In order to provide the distribution of the concentration process $C$ by using that $C_. = |X_.|^{\gamma + 1}$, the following lemma provides first the covariance function of the fractional Ornstein-Uhlenbeck process $X$.
%

% Lemma : Auxiliary distributions.

%
\begin{lemma}\label{auxiliary_distributions}
$B^H(\vartheta)$ is a centered Gaussian process of covariance function $R_{H,\vartheta}$ such that :
\begin{displaymath}
R_{H,\vartheta}(s,t) =
\alpha_H(1-\beta)^2
\int_{0}^{s}\int_{0}^{t}
|u - v|^{2(H - 1)}e^{\upsilon(1-\beta)(u + v)}dudv
\end{displaymath}
for every $s,t\in\mathbb R_+$. Then, the covariance function $R_X$ of the fractional Ornstein-Uhlenbeck process $X$ satisfies :
\begin{displaymath}
R_X(s,t) =
\alpha_H\sigma^2(1-\beta)^2
\int_{0}^{s}\int_{0}^{t}
|u - v|^{2(H - 1)}e^{-\upsilon(1-\beta)[(t - u) + (s - v)]}dudv
\end{displaymath}
for every $s,t\in\mathbb R_+$.
\end{lemma}
\noindent
The following proposition provides the finite-dimensional distributions of the concentration process $C$.
%

% Proposition : Distribution of the concentration process.

%
\begin{proposition}\label{distribution_concentration_process}
For every $n\in\mathbb N^*$ and $t_1,\dots,t_n\in\mathbb R_+$, the distribution of the random vector $(C_{t_1},\dots,C_{t_n})$ admits a density $\chi_n$ with respect to the Lebesgue measure on $(\mathbb R^n,\mathcal B(\mathbb R^n))$ such that :
\begin{eqnarray*}
 \chi_n(x_1,\dots,x_n) & = &
 \frac{2^n(1-\beta)^n\mathbf 1_{\mathbb R_{+}^{n}}(x_1,\dots,x_n)}{(2\pi)^{n/2}\sqrt{|\det(R_n)|}}
 \prod_{i = 1}^{n}x_{i}^{-\beta}\times\\
 & &
 \exp\left[-\frac{1}{2}
 \left[
 R_{n}^{-1}
 \left[
 \begin{pmatrix}
 x_{1}^{1-\beta}\\
 \vdots\\
 x_{n}^{1-\beta}
 \end{pmatrix}
 -V_n\right]\right]
 \cdot
 \left[
 \begin{pmatrix}
 x_{1}^{1-\beta}\\
 \vdots\\
 x_{n}^{1-\beta}
 \end{pmatrix}
 -V_n
 \right]\right]\\
 & &
 \textrm{$;$ }
 (x_1,\dots,x_n)\in\mathbb R^n
\end{eqnarray*}
where, $R_n\in S_n(\mathbb R_+)$ and $V_n\in\mathbb R^n$ satisfy
\begin{displaymath}
R_n(i,j) :=
\alpha_H\sigma^2(1-\beta)^2
\int_{0}^{t_i}\int_{0}^{t_j}
|u - v|^{2(H - 1)}e^{-\upsilon(1-\beta)\left[(t_j - u) + (t_i - v)\right]}dudv
\end{displaymath}
and
\begin{displaymath}
V_n(i) :=
C_{0}^{1-\beta}
e^{-\upsilon(1-\beta)t_i}
\end{displaymath}
for every $i,j\in\{1,\dots,n\}$.
\end{proposition}
\noindent
It is a straightforward application of Lemma \ref{auxiliary_distributions} together with N. Marie \cite{MARIE12}, Proposition 5.1.
\\
\\
The following proposition allows to control, in probability, the uniform distance between the process $X$ and the solution $X^{\det}$ of the associated ordinary differential equation
\begin{displaymath}
X_{t}^{\det} =
X_0 -
\upsilon(1-\beta)\int_{0}^{t}X_{s}^{\det}ds
\textrm{ ; }
t\in\mathbb R_+.
\end{displaymath}
%

% Proposition : Upper-bound for choosing the variance.

%
\begin{proposition}\label{variance_upper_bound}
For every $x > 0$ and $T > 0$,
\begin{displaymath}
\mathbb P\left(\left\|X - X^{\det}\right\|_{\infty,T} > x\right)
\leqslant
2\exp\left[
-\frac{x^2}{2\sigma^2R_{H,\vartheta}(T,T)}\right].
\end{displaymath}
\end{proposition}
\noindent
For a level $\lambda\in ]0,1[$, in order to ensure with probability greater than $1-\lambda$ that $|X_t - X_{t}^{\det}|\leqslant x\in\mathbb R_{+}^{*}$ for every $t\in [0,T]$, it is sufficient to assume that $\sigma^2\in [0,M(\lambda,x,H)]$ with
\begin{displaymath}
M(\lambda,x,H) :=
\frac{x^2}{2R_{H,\vartheta}(T,T)\log(2/\lambda)}
\end{displaymath}
by Proposition \ref{variance_upper_bound}.
%

% Corollary : Consequence on the concentration process.

%
\begin{corollary}\label{concentration_process_consequence}
For every $x > 0$ and $T > 0$,
\begin{displaymath}
\mathbb P\left[
\forall t\in [0,T]
\textrm{$,$ }
C_t\in\left[0,2^{\gamma}(C_{t}^{\det} + x^{\gamma + 1})\right]\right]
\geqslant
1 - 2\exp\left[-\frac{x^2}{2\sigma^2R_{H,\vartheta}(T,T)}\right]
\end{displaymath}
where $C_{t}^{\det} := |X_{t}^{\det}|^{\gamma + 1}\textrm{ $;$ }t\in\mathbb R_+$.
\end{corollary}
\noindent
An illustration of how one can use Proposition \ref{variance_upper_bound} and Corollary \ref{concentration_process_consequence} in practice is provided at Subsection 3.1.
%

% Subsection : Estimator of the elimination constant and ergodic theorem.

%
\subsection{Ergodic theorem and estimator of the elimination constant}
The following proposition is an extension of existing ergodic theorems for the fractional Ornstein-Uhlenbeck process $X$. Let $Y$ be the stochastic process defined by :
\begin{displaymath}
Y_t :=
\sigma(1-\beta)
\int_{-\infty}^{t}e^{-\upsilon(1-\beta)(t - s)}dB_{s}^{H}
\textrm{ ; }
t\in\mathbb R_+.
\end{displaymath}
%

% Proposition : An ergodic theorem for the fractional Ornstein-Uhlenbeck process.

%
\begin{proposition}\label{ergodic_theorem_OU}
Let $f :\mathbb R\rightarrow\mathbb R$ be a continuous function such that :
\begin{eqnarray*}
 \exists n\in\mathbb N^*
 \textrm{$,$ }
 \exists (a_1,b_1,c_1),\dots,(a_n,b_n,c_n)\in\mathbb R_{+}^{*}\times\mathbb R_{+}^{2} :
 \forall x,\varepsilon & \in & \mathbb R\textrm{$,$ }\\
 |f(x + \varepsilon) - f(x)| & \leqslant &
 \sum_{i = 1}^{n}
 c_i(1 + |x|)^{b_i}|\varepsilon|^{a_i}.
\end{eqnarray*}
Then,
\begin{displaymath}
\frac{1}{T}
\int_{0}^{T}
f(X_t)dt
\xrightarrow[T\rightarrow\infty]{\textrm{a.s.}}
\mathbb E\left[f(Y_0)\right] < \infty.
\end{displaymath}
\end{proposition}
\noindent
\noindent
With notations of Proposition \ref{ergodic_theorem_OU}, put $f(x) := x^n$ ; $x\in\mathbb R$. For every $x,\varepsilon\in\mathbb R$,
\begin{eqnarray}
 |f(x +\varepsilon) - f(x)| & \leqslant &
 \sum_{i = 0}^{n - 1}
 \begin{pmatrix}
 n\\
 i
 \end{pmatrix}
 |x|^i|\varepsilon|^{n - i}
 \nonumber\\
 \label{power_ergodic_theorem_OU}
 & \leqslant &
 \sum_{i = 0}^{n - 1}
 p_i(x)|\varepsilon|^{a_i}
\end{eqnarray}
where $a_i := n - i$ and
\begin{displaymath}
p_i(x) :=
\begin{pmatrix}
n\\
i
\end{pmatrix}
(1 + |x|)^i
\textrm{ ; }
x\in\mathbb R
\textrm{, }
i = 0,\dots,n - 1.
\end{displaymath}
Then, by Proposition \ref{ergodic_theorem_OU} :
\begin{eqnarray}
 \lim_{T\rightarrow\infty}
 \frac{1}{T}
 \int_{0}^{T}X_{t}^{n}dt & = &
 \mathbb E\left(Y_{0}^{n}\right)
 \textrm{ $\mathbb P$-a.s.}
 \nonumber\\
 \label{ergodic_Hu_Nualart}
 & = &
 \left\{
 \begin{array}{rcl}
 \displaystyle{
 \frac{n!\sigma^{n}(1-\beta)^{n - nH}\upsilon^{-nH}H^{n/2}\Gamma^{n/2}(2H)}{2^{n/2}(n/2)!}}
 & \textrm{if} & n\in 2\mathbb N^*\\
 0 & \textrm{if} & n\in\mathbb N^*-(2\mathbb N^*)
 \end{array}\right.
\end{eqnarray}
by Y. Hu and D. Nualart \cite{HN10}, Lemma 5.1. For $n = 2$, (\ref{ergodic_Hu_Nualart}) coincides with Y. Hu and D. Nualart \cite{HN10}, Lemma 3.3.
\\
\\
Assume that values of parameters $H$ and $\sigma$ are known. For $T > 0$ arbitrarily chosen, consider
\begin{displaymath}
\widehat{\upsilon}_T :=
\frac{1}{1-\beta}
\left[
\frac{1}{\sigma^2(1-\beta)^2H\Gamma(2H)T}
\int_{0}^{T}X_{t}^{2}dt
\right]^{-1/(2H)}.
\end{displaymath}
%

% Proposition : Strong estimator of $\upsilon$.

%
\begin{proposition}\label{strong_upsilon_estimator}
$\widehat{\upsilon}_T$ is a strongly consistent estimator of $\upsilon$.
\end{proposition}
\noindent
It is a straightforward consequence of (\ref{ergodic_Hu_Nualart}) for $n = 2$. The estimator $\widehat{\upsilon}_T$ was studied by Y. Hu and D. Nualart at \cite{HN10}, Section 4. They completed Proposition \ref{strong_upsilon_estimator} by a \textit{central limit theorem} when $H\in ]1/2,3/4[$ (cf. \cite{HN10}, Theorem 4.1).
%

% Subsection : Estimators of the Hurst parameter and of the volatility constant.

%
\subsection{Estimators of the Hurst parameter and of the volatility constant}
Assume that the concentration process $C$ is discretely observed at times $t_0,\dots,t_n$, where $n\in\mathbb N^*$, $t_k := k\delta_n$ for every $k\in\{0,\dots,n\}$, and $(\delta_n,n\in\mathbb N)$ is a $\mathbb R_{+}^{*}$-valued sequence such that
\begin{displaymath}
\lim_{n\rightarrow\infty}
\delta_n = 0
\textrm{ and }
\lim_{n\rightarrow\infty}
n\delta_n =\infty.
\end{displaymath}
Proposition \ref{H_estimator} provides a strongly consistent estimator, easy to implement, of the Hurst parameter $H$ coming from J. Istas and G. Lang \cite{IL97}. Proposition \ref{sigma_estimator} provides an associated consistent estimator of $\sigma$.
\\
Proposition \ref{weak_upsilon_estimator} provides a weakly consistent estimator of $\upsilon$ for unknown values of parameters $H$ and $\sigma$.
%

% Proposition : Estimator of $H$.

%
\begin{proposition}\label{H_estimator}
Consider
\begin{displaymath}
\widehat{H}_n :=
\frac{1}{2}
\log_2\left(
\frac{\displaystyle{\sum_{k = 2}^{n - 2}\left|X_{t_{k + 2}} - 2X_{t_k} + X_{t_{k - 2}}\right|^2}}
{\displaystyle{\sum_{k = 1}^{n - 1}\left|X_{t_{k + 1}} - 2X_{t_k} + X_{t_{k - 1}}\right|^2}}\right).
\end{displaymath}
$\widehat{H}_n$ is a strongly consistent estimator of $H$.
\end{proposition}
%

% Proposition : Estimator $\sigma$.

%
\begin{proposition}\label{sigma_estimator}
Consider $a_0 := -1/4$, $a_1 := 1/2$, $a_2 := -1/4$ and
\begin{displaymath}
\widehat{\sigma}_n :=
\frac{1}{1-\beta}
\left(-\frac{1}{8}\times
\frac{\displaystyle{\sum_{k = 1}^{n - 1}\left|X_{t_{k + 1}} - 2X_{t_k} + X_{t_{k - 1}}\right|^2}}
{\displaystyle{\sum_{k,l = 0}^{2}
a_ka_l|k - l|^{2\widehat{H}_n}\delta_{n}^{2\widehat{H}_n}}}
\right)^{1/2}.
\end{displaymath}
$\widehat{\sigma}_n$ is a strongly consistent estimator of $\sigma$.
\end{proposition}
\noindent
Refer to A. Brouste and S. Iacus \cite{BI13}, Theorem 1, based on J. Istas and G. Lang \cite{IL97}, Theorem 3, for a proof of propositions \ref{H_estimator} and \ref{sigma_estimator}.
\\
The R-package Yuima, developed by A. Brouste and S. Iacus, allows to compute estimations of $(H,\sigma)$ via $(\widehat{H}_n,\widehat{\sigma}_n)$.
%

% Proposition : Weak estimator of $\upsilon$.

%
\begin{proposition}\label{weak_upsilon_estimator}
Consider
\begin{displaymath}
\widehat{\upsilon}_{n}^{*} :=
\frac{1}{1-\beta}\left[
\frac{1}{\widehat{\sigma}_{n}^{2}(1-\beta)^2\widehat H_n\Gamma(2\widehat{H}_n)n}
\sum_{k = 0}^{n - 1}X_{k\delta_n}^{2}
\right]^{-1/(2\widehat{H}_n)}.
\end{displaymath}
$\widehat{\upsilon}_{n}^{*}$ is a weakly consistent estimator of $\upsilon$.
\end{proposition}
%

% Section : Numerical simulations and pharmacokinetics.

%
\section{Numerical simulations and pharmacokinetics}
\noindent
For small sets of observations, the first subsection provides a qualitative procedure for choosing parameters $H$, $\sigma$ and $\beta$. Proposition \ref{variance_upper_bound} is the cornerstone of the procedure. The second subsection illustrates the convergence of estimators provided at Section 2. The relationship between the estimations quality and the size/length of the sample is discussed.
%

% Subsection : A qualitative procedure for choosing $H$, $\sigma$ and $\beta$.

%
\subsection{A qualitative procedure for choosing $H$, $\sigma$ and $\beta$}
Consider $n\in\mathbb N^*$ and $(t_1,\dots,t_n)\in\mathbb R_{+}^{n}$ satisfying $t_1 <\dots < t_n\leqslant T$. Throughout this subsection, assume that concentrations have been observed at times $t_1,\dots,t_n$. These concentrations $c_1,\dots,c_n$ provide observations $x_1,\dots,x_n$ of the fractional Ornstein-Uhlenbeck process $X$ by putting $x_i = c_{i}^{1-\beta}$ ; $i = 1,\dots,n$.
\\
\\
Consider the following values of the other parameters involved in equation (\ref{main_model}), coming from Y. Jacomet \cite{JACOMET89}, Chapitre II.3 :
\begin{center}
\begin{tabular}{|l| l |}
\hline
Parameters & Values\\
\hline
\hline
$T$ & 3h\\
$\upsilon$ & 3.5$\textrm{h}^{-1}$\\
$C_0$ & 1g\\
\hline
$n$ & 500\\
\hline
\end{tabular}
\end{center}
In order to choose $H$, $\sigma$ and $\beta$, a qualitative procedure is provided by using these values as an example. That method is simple and doesn't require a lot of observations of the concentration process.
\\
\\
On one hand, as mentioned at Proposition \ref{variance_upper_bound}, for a level $\lambda\in ]0,1[$, in order to ensure with probability greater than $1-\lambda$ that $|X_t - X_{t}^{\det}|\leqslant x\in\mathbb R_{+}^{*}$ for every $t\in [0,T]$, it is sufficient to assume that $\sigma^2\in [0,M(\lambda,x,H)]$ with
\begin{displaymath}
M(\lambda,x,H) :=
\frac{x^2}{2R_{H,\vartheta}(T,T)\log(2/\lambda)}.
\end{displaymath} 
Moreover, by Corollary \ref{concentration_process_consequence} :
\begin{displaymath}
\mathbb P\left[
\forall t\in [0,T]\textrm{, }
C_t\in\left[0,2^{\gamma}\left(C_{t}^{\det} + x^{\gamma + 1}\right)\right]
\right]
\geqslant 1-\lambda.
\end{displaymath}
On the other hand, as mentioned in introduction, the H\"older regularity of the concentration process paths is continuously controlled by the Hurst parameter $H$.
\\
\\
For $H = 0.9$ and $\beta = 0.9$, the following array provides the values of $M(\lambda,x,H)$ for usual levels $\lambda = 0.01,0.05,0.10$ :
\begin{center}
\begin{tabular}{|l| l l l |}
\hline
$x^{\gamma + 1}$ | $\lambda$ & 0.01 & 0.05 & 0.10\\
\hline
0.1 & 0.26 & 0.38 & 0.46\\
0.2 & 0.30 & 0.43 & 0.53\\
0.4 & 0.36 & 0.50 & 0.61\\
\hline
\end{tabular}
\end{center}
On the two following figures, the paths of the process $C$ are respectively plotted for extreme cases $\sigma^2 = 0.26$ and $\sigma^2 = 4 > 0.61$. The concentration process paths are plotted in black and the associated deterministic model is plotted in red :
\begin{figure}[htbp]
\begin{minipage}[c]{.45\linewidth}
\begin{center}
\includegraphics[scale = 0.25]{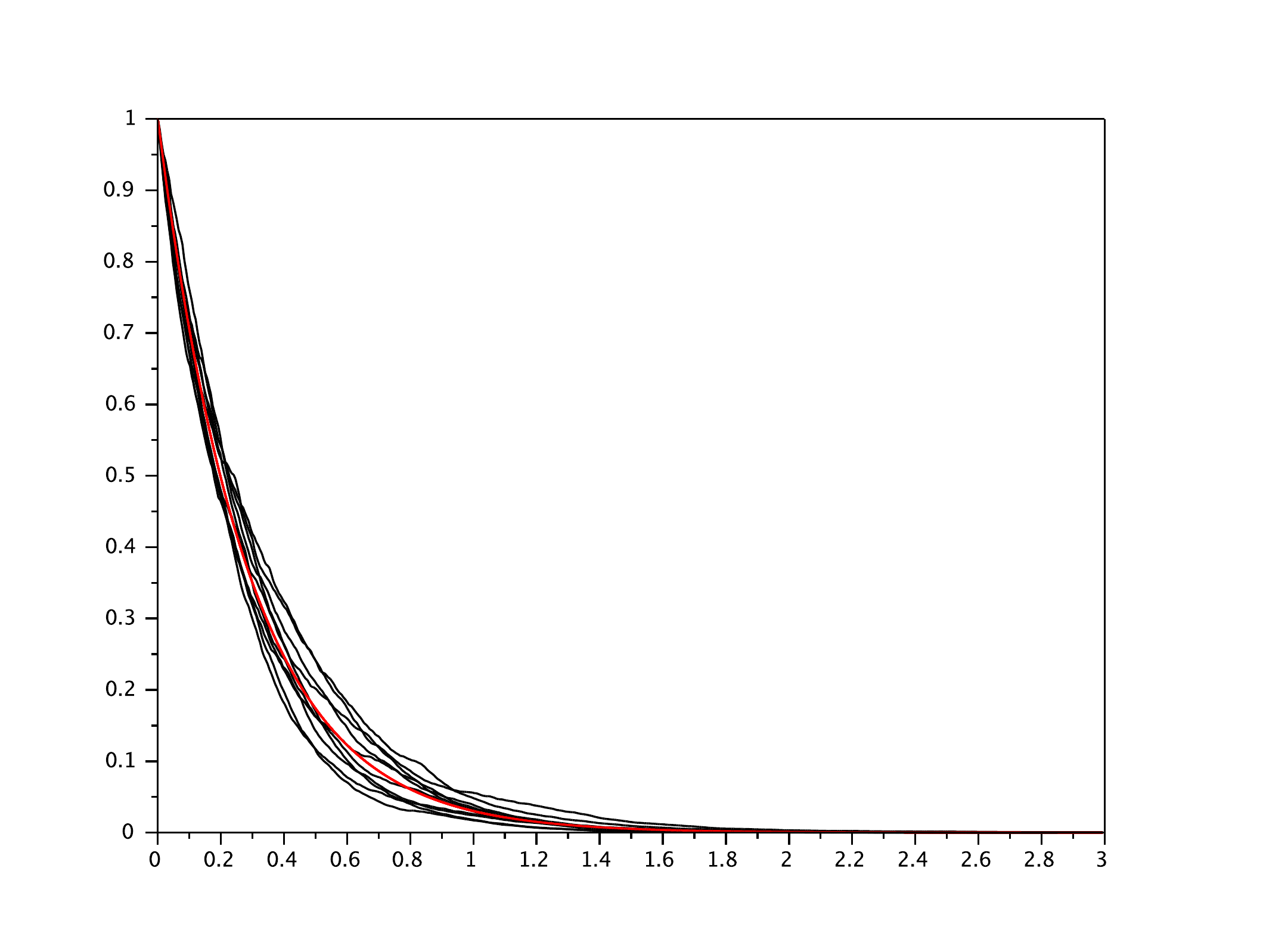}
\caption{$H = 0.9$ and $\sigma^2 = 0.26$}
\label{fig:image1}
\end{center}
\end{minipage}
\hfill
\begin{minipage}[c]{.45\linewidth}
\begin{center}
\includegraphics[scale = 0.25]{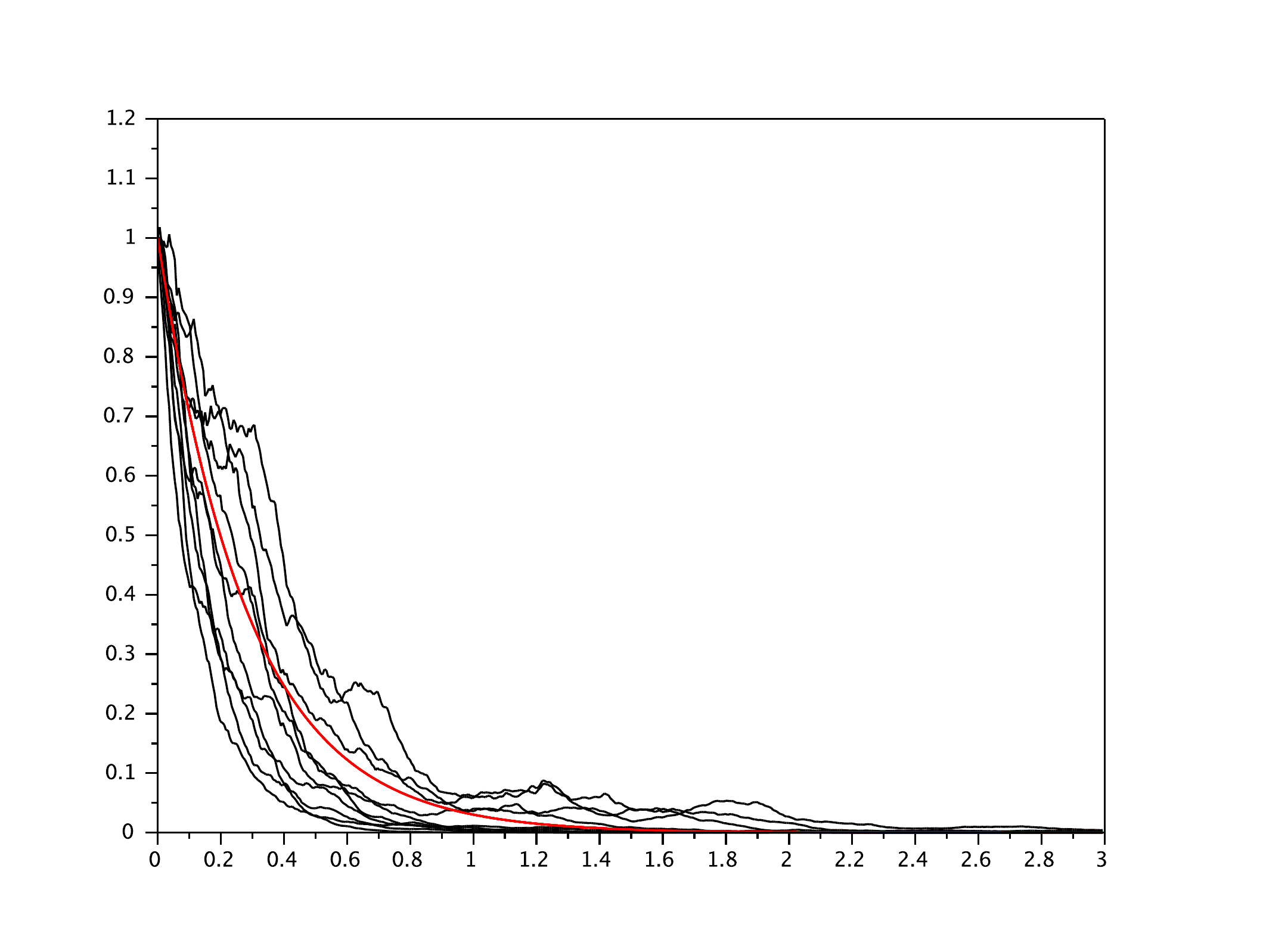}
\caption{$H = 0.9$ and $\sigma^2 = 4 > 0.61$}
\label{fig:image2}
\end{center}
\end{minipage}
\end{figure}
\newline
For $H = 0.6$ and $\beta = 0.9$, the following array provides the values of $M(\lambda,x,H)$ for usual levels $\lambda = 0.01,0.05,0.10$ :
\begin{center}
\begin{tabular}{|l| l l l |}
\hline
$x^{\gamma + 1}$ | $\lambda$ & 0.01 & 0.05 & 0.10\\
\hline
0.1 & 0.70 & 1.00 & 1.23\\
0.2 & 0.80 & 1.15 & 1.42\\
0.4 & 0.92 & 1.32 & 1.63\\
\hline
\end{tabular}
\end{center}
On the two following figures, paths of the process $C$ are plotted as for \mbox{$H = 0.9$ :}
\begin{figure}[htbp]
\begin{minipage}[c]{.45\linewidth}
\begin{center}
\includegraphics[scale = 0.25]{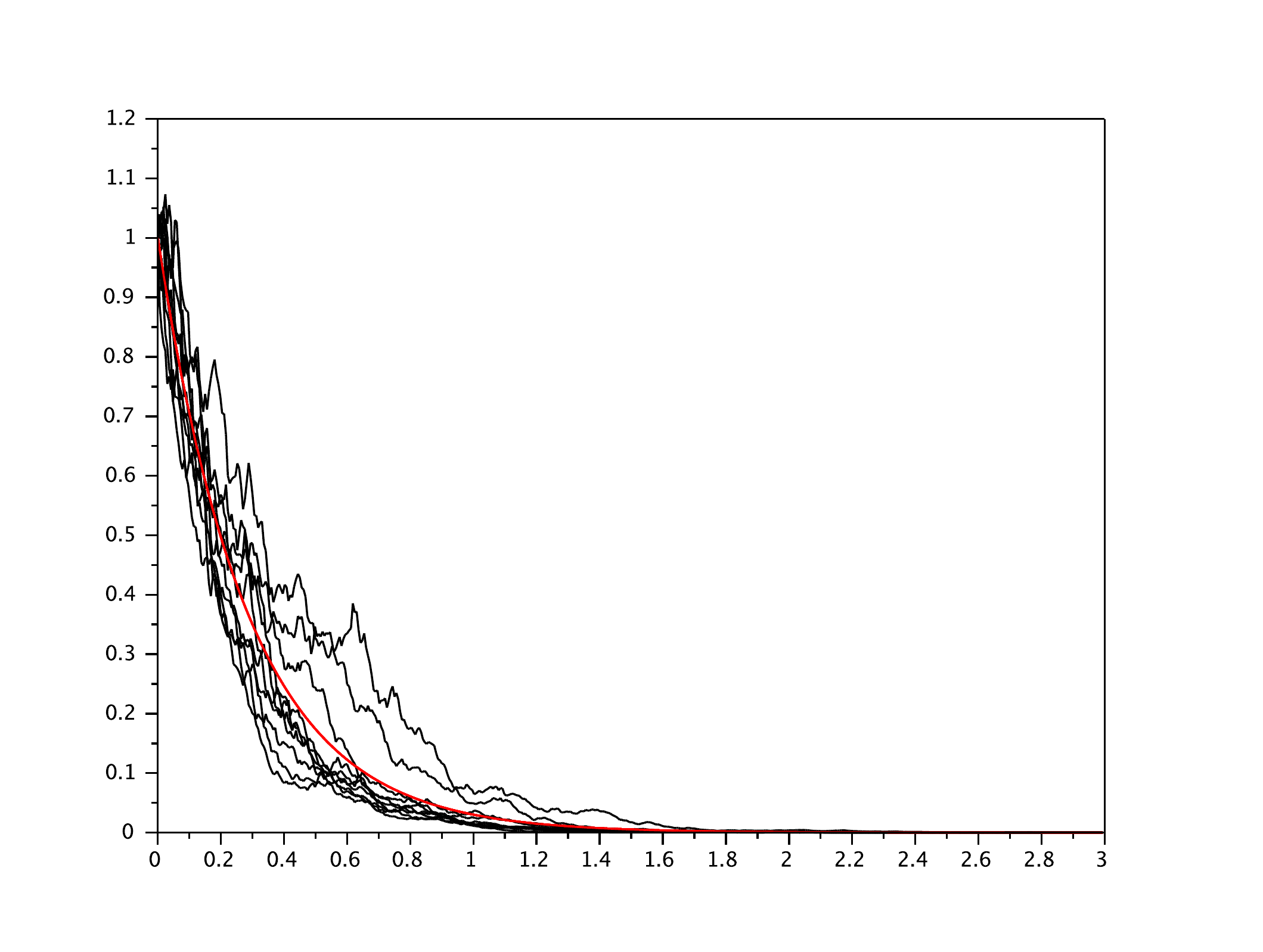}
\caption{$H = 0.6$ and $\sigma^2 = 0.70$}
\label{fig:image1}
\end{center}
\end{minipage}
\hfill
\begin{minipage}[c]{.45\linewidth}
\begin{center}
\includegraphics[scale = 0.25]{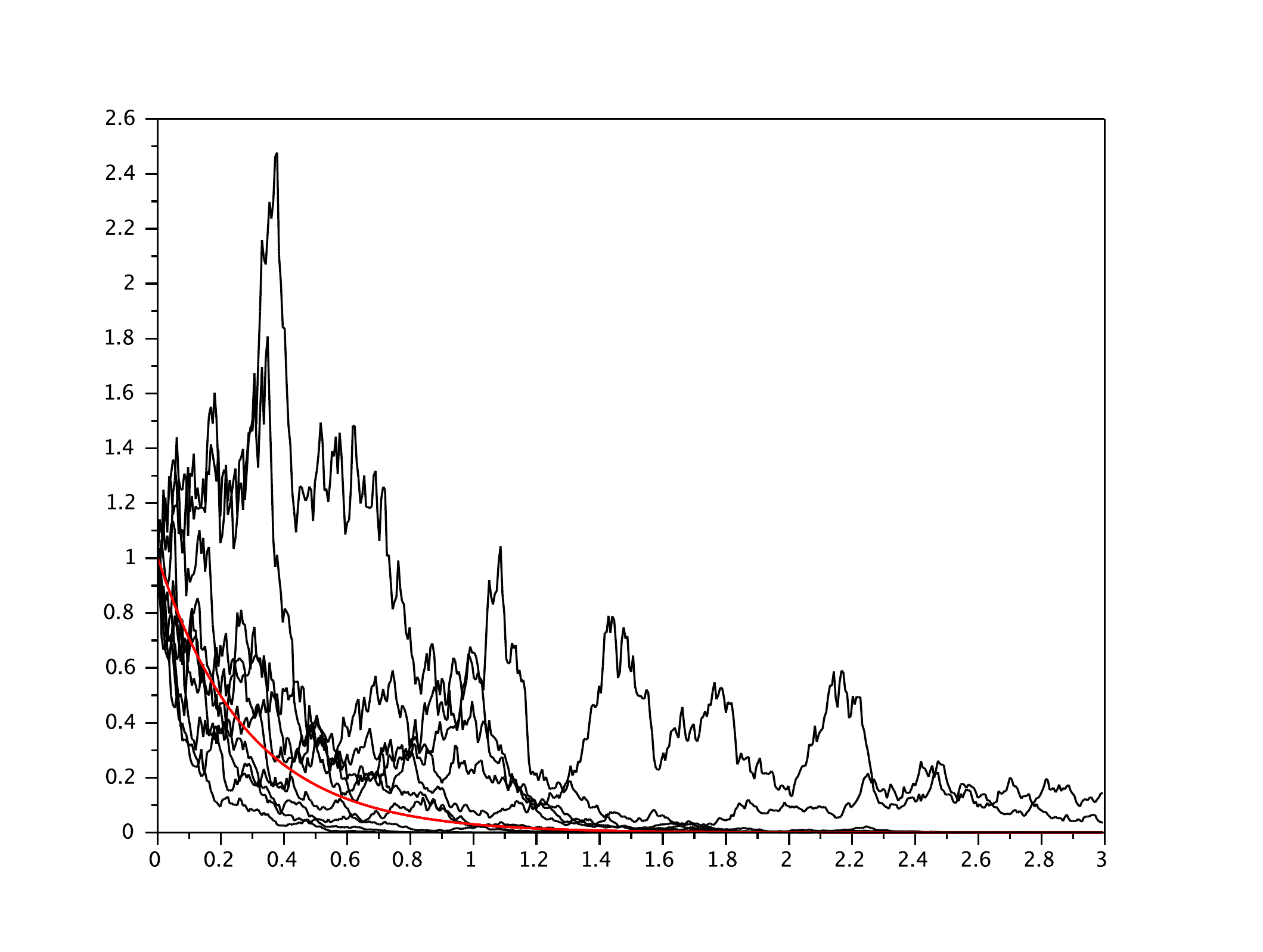}
\caption{$H = 0.6$ and $\sigma^2 = 4 > 1.63$}
\label{fig:image2}
\end{center}
\end{minipage}
\end{figure}
\newpage
\noindent
On one hand, in order to model the concentration process realistically, one should take $H = 0.9$ and, for instance :
\begin{eqnarray*}
 \sigma^2 & \in &
 ]0;M(0.01,0.2^{1-\beta},0.9)]\\
 & = &
 ]0,0.30].
\end{eqnarray*}
Indeed,
\begin{itemize}
 \item For $H = 0.9$ with $\sigma^2 = 4 > M(0.10,0.4^{1-\beta},0.9) = 0.61$, the concentration process paths seem locally regular enough, but not globally.
 \item For $H = 0.6$ with $\sigma^2 = M(0.01,0.1^{1-\beta},0.6) = 0.70$, the concentration process paths seem globally regular enough, but not locally.
 \item For $H = 0.6$ with $\sigma^2 = 4 > M(0.10,0.4^{1-\beta},0.6) = 1.63$, the concentration process paths seem not regular enough locally and globally.
\end{itemize}
Also, $\beta = 0.9$ seems to be a good choice. Indeed, if $\beta\in [0,0.8]$, for every
\begin{displaymath}
\sigma\in
\left[M\left(0.01,0.1^{1-\beta},0.9\right) ; M\left(0.10,0.4^{1-\beta},0.6\right)\right],
\end{displaymath}
the concentration process paths seem not significantly perturbed with respect to the associated deterministic model. Then, to take $\beta = 0.9$ ensures that the value of the parameter $\sigma$ can be chosen such that the following realistic condition is \mbox{satisfied :}
\begin{displaymath}
\mathbb P\left(\forall t\in [0,T]
\textrm{, }X_t\in\left[
X_{t}^{\det} - 0.2^{1-\beta},X_{t}^{\det} + 0.2^{1-\beta}\right]\right)
\geqslant 0.99.
\end{displaymath}
On the observed concentrations $c_1,\dots, c_n$, the following procedure allows to choose $H$, $\sigma$ and $\beta$ qualitatively :
\begin{itemize}
 \item\textbf{Step 1.} Take $H\in ]0.5,1[$ sufficiently close to $1$ as $0.9$.
 \item\textbf{Step 2.} Take $\beta\in ]0,1[$.
 \item\textbf{Step 3.} Choose a standard level $\lambda\in ]0,1[$ as $0.01$ or $0.05$, and take for instance
 \begin{eqnarray*}
  x & := &
  \max_{i = 1}^{n}
  \left|c_{i}^{1-\beta} -
  (C_0e^{-\upsilon t_i})^{1-\beta}\right|\\
  & = &
  \max_{i = 1}^{n}
  \left|x_i -
  C_{0}^{1-\beta}e^{-\upsilon(1-\beta)t_i}\right|.
 \end{eqnarray*}
 Then, compute $M(\lambda,x,H)$.
 \\
 If the value of $\upsilon$ is unknown, since paths of the concentration process have to be moderately perturbed with respect to the associated deterministic model, it can be approximated by linear regression as in Y. Jacomet \cite{JACOMET89} (see Subsection 3.2).
 \item\textbf{Step 4.} Take $\sigma^2\in ]0;M(\lambda,x,H)]$ such that the concentration process paths seem regular enough locally and globally to model the elimination of the administered drug.
 \\
 If the concentration process paths are not significantly perturbed with respect to the associated deterministic model for usual levels $\lambda\in ]0,1[$, then go to the second step and choose a greater value of the parameter $\beta$. If the concentration process paths are not globally regular enough for standard levels $\lambda\in ]0,1[$, then go to the second step and choose a smaller value of the parameter $\beta$.
\end{itemize}
%

% Subsection : Parameters estimations.

%
\subsection{Parameters estimation}
Throughout this subsection, assume that the concentration process $C$ has been discretely observed at times $t_0,\dots,t_n$, where $n\in\mathbb N^*$, $t_k := k\delta_n$ for every $k\in\{0,\dots,n\}$, and $(\delta_n,n\in\mathbb N)$ is a $\mathbb R_{+}^{*}$-valued sequence such that
\begin{displaymath}
\lim_{n\rightarrow\infty}
\delta_n = 0
\textrm{ and }
\lim_{n\rightarrow\infty}
n\delta_n =\infty.
\end{displaymath}
Consider the following values of parameters involved in equation (\ref{main_model}) :
\begin{center}
\begin{tabular}{|l| l |}
\hline
Parameters & Values\\
\hline
\hline
$T$ & $n\delta_n$ ; $n = 10,\dots,10^3$\\
$\beta$ & 0\\
$\upsilon$ & 1.5$\textrm{h}^{-1}$\\
$H$ & 0.9\\
$\sigma^2$ & $0.26$\\
$C_0$ & 1g\\
\hline
\end{tabular}
\end{center}
\noindent
The two following figures illustrate the convergence of estimators $\widehat{\upsilon}_n$ and $\widehat{H}_n$ provided at propositions \ref{strong_upsilon_estimator} and \ref{H_estimator} respectively. For every $n$ belonging to $\{10,\dots,10^3\}$, the concentration process $C$ is simulated at times $t_0,\dots,t_n$ and estimators $\widehat\upsilon_n$ and $\widehat H_n$ are computed with these simulated observations denoted by $c_1,\dots,c_n$. Estimations are plotted in black and parameters values are plotted in red :
\begin{figure}[htbp]
\begin{minipage}[c]{.45\linewidth}
\begin{center}
\includegraphics[scale = 0.30]{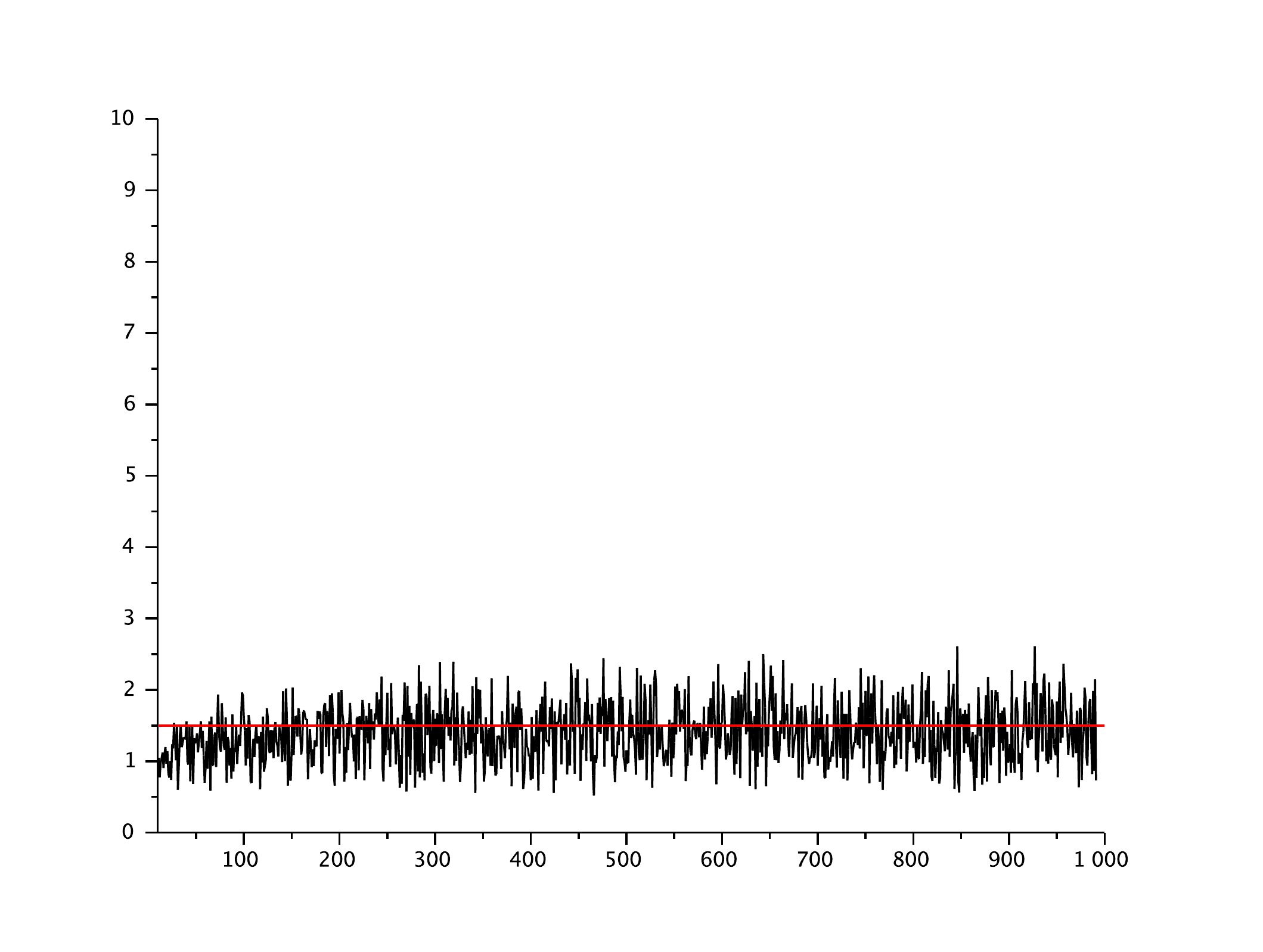}
\caption{$\widehat{\upsilon}_n$ for $n = 10,\dots,10^3$}
\label{fig:image1}
\end{center}
\end{minipage}
\hfill
\begin{minipage}[c]{.45\linewidth}
\begin{center}
\includegraphics[scale = 0.30]{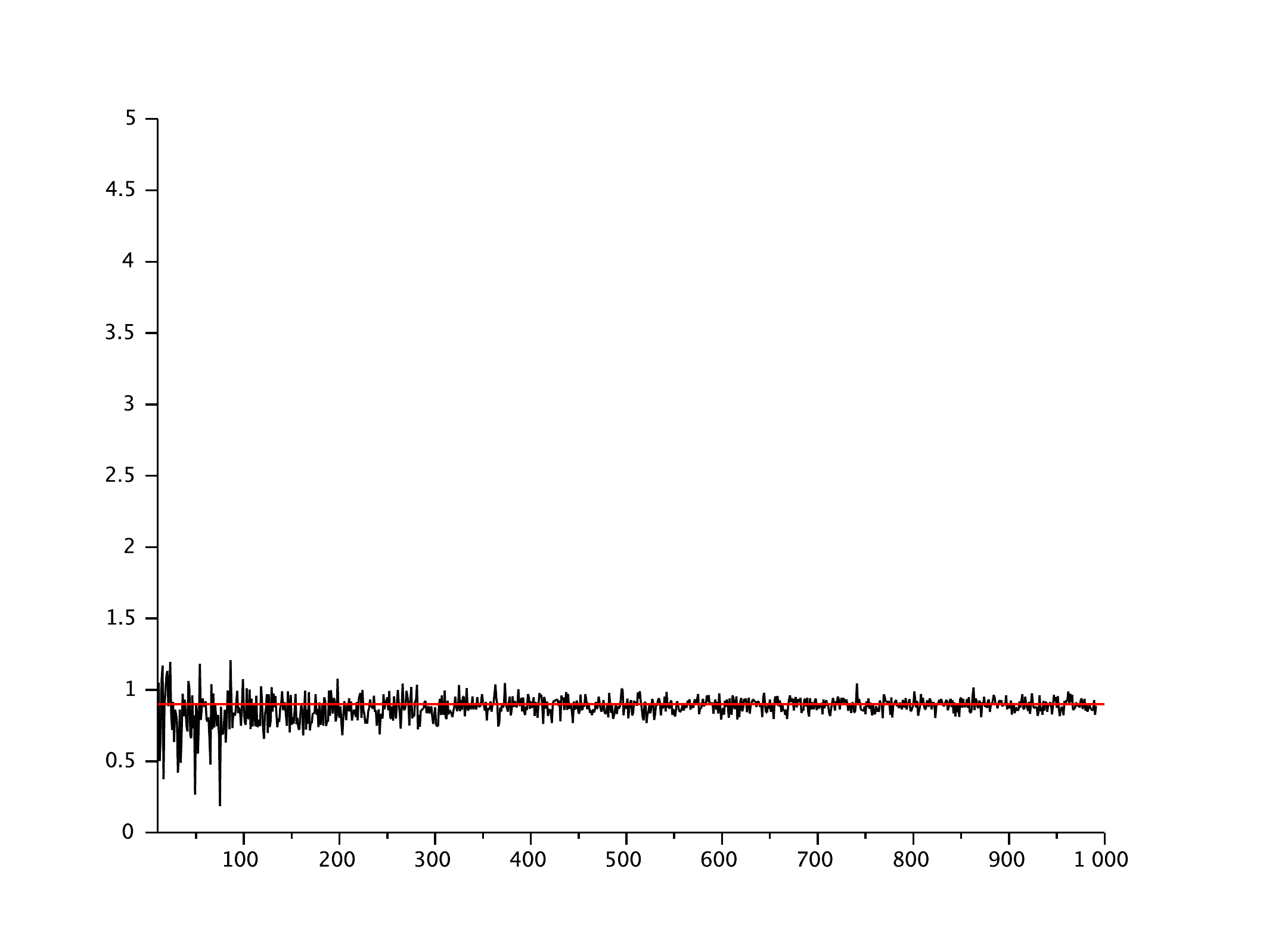}
\caption{$\widehat H_n$ for $n = 10,\dots,10^3$}
\label{fig:image2}
\end{center}
\end{minipage}
\end{figure}
\newline
The estimator $\widehat{\upsilon}_n$ converges slowly to the elimination constant $\upsilon$. Then, if the number $n$ of observations is insufficient, since paths of the concentration process have to be moderately perturbed with respect to the associated deterministic model, one can take as in Y. Jacomet \cite{JACOMET89} :
\begin{displaymath}
\upsilon
\approx
-\frac{1}{1-\beta}\times\frac{\textrm{cov}\left[
t_1,\dots,t_n ; \log(c_1),\dots,\log(c_n)\right]}{\textrm{var}(t_1,\dots,t_n)}.
\end{displaymath}
The qualitative procedure provided at Subsection 3.1 is also an alternative for choosing $H$ and $\sigma$ with few observed concentrations.
%

% Section : Discussion and perspectives.

%
\section{Discussion and perspectives}
\noindent
The stochastic model studied in this paper is a natural extension of usual deterministic models used in pharmacokinetics, it has smooth enough paths to take realistically in account the random component of the elimination process, and its explicit expression together with Decreusefond-Lavaud method allow to simulate it easily. As mentioned at Section 3, estimators of parameters $\upsilon$, $H$ and $\sigma$ provide good estimations for large sets of observed concentrations. For small sets of observations, the qualitative procedure described at Subsection 3.1 is simple and seems quite efficient. For these reasons, the model could be used in clinical applications.
\\
\\
Assume that the therapeutic response $R_t$ to the administered drug at time $t\in [0,\tau_0]$ satisfies $R_t := F(C_t,O_t)$, where $F\in C^1(\mathbb R^2;\mathbb R)$ and $O$ is a stochastic process with $\mathbb R$-valued paths that doesn't depend on the initial concentration $C_0 = A/V$. The random variable $C_t$ is derivable with respect to $C_0 > 0$ and
\begin{displaymath}
\frac{\partial R_t}{\partial C_0} =
C_{0}^{-\beta}\left[C_{0}^{1-\beta} + \sigma B_{t}^{H}(\vartheta)\right]^{\gamma}e^{-\upsilon t}
\partial_1F(C_t,O_t).
\end{displaymath}
Differential calculus arguments could then allow to compute the dose that maximises the therapeutic response $R_t$ for some well chosen functions $F$ and well chosen stochastic processes $O$.
\\
\\
Since the stochastic process $C$ seems to model the elimination process more realistically than the deterministic function $C^{\det}$, the perspective of clinical applications described above could be interesting for potentially toxic drugs.
\\
For instance, the elimination of the ketamine, that can be neurotoxic but more effective than classic antidepressant in the treatment of major depressive disorders (cf. G.E. Correll and G.E. Futter \cite{CF06}), could be modeled by the stochastic process studied in that paper. To choose $F$ and $O$ such that $R_t$ models the Hamilton rating scale or the Beck depression inventory at time $t$ could allow to compute the dose of ketamine maximizing its antidepressant effect and minimizing its neurotoxic effect.
\appendix
%

% Appendix A : Fractional Brownian motion.

%
\section{Fractional Brownian motion}
\noindent
Essentially inspired by D. Nualart \cite{NUALART06} and, L. Decreusefond and A.S. Ust\"unel \cite{DU99}, this section gives basics on the fractional Brownian motion $B^H$ of Hurst parameter $H\in ]1/2,1[$, its reproducing kernel Hilbert space and the fractional Young/Wiener integral with respect to $B^H$.
\\
On Gaussian processes, the reader can refer to J. Neveu \cite{NEVEU68}.
\\
\\
For a time $T > 0$ arbitrarily chosen, consider
\begin{displaymath}
\Delta_T :=
\left\{
(s,t)\in [0,T]^2 :
s < t\right\}.
\end{displaymath}
%

% Definition : The fractional Brownian motion.

%
\begin{definition}\label{fBm}
A fractional Brownian motion of Hurst parameter $H\in ]0,1]$ is a centered Gaussian process $B^H$ of covariance function $R_H$ defined by :
\begin{displaymath}
R_H(s,t) :=
\frac{1}{2}(
s^{2H} + t^{2H} - |t - s|^{2H})
\textrm{ $;$ }
s,t\in\mathbb [0,T].
\end{displaymath}
\end{definition}
\noindent
The process $B^H$ is a semi-martingale if and only if $H = 1/2$ (cf. \cite{NUALART06}, Proposition 5.1.1). Then, it is not possible to integrate with respect to $B^H$ in the sense of It\^o. However, since
\begin{displaymath}
\mathbb E\left(|B_{t}^{H} - B_{s}^{H}|^2\right) = |t - s|^{2H}
\end{displaymath}
for every $s,t\in [0,T]$, the Kolmogorov continuity criterion ensures that $B^H$ has $\alpha$-H\"older continuous paths with $\alpha\in ]0,H[$. Therefore, for any stochastic process $X$ with $\beta$-H\"older continuous paths such that $\alpha +\beta > 1$, it is possible to integrate $X$ with respect to $B^H$ in the sense of Young.
\\
About the Young integral, that extends the classic Riemann-Stielj\`es integral, the reader can refer to A. Lejay \cite{LEJAY10}.
\\
\\
In the sequel, assume that $H\in ]1/2,1[$ and put $\alpha_H := H(2H - 1)$. The vector space
\begin{displaymath}
\mathcal H :=
\left\{
h\in L^2([0,T];dt) :
\alpha_H
\int_{0}^{T}
\int_{0}^{T}
|t - s|^{2(H - 1)}h(s)h(t)dsdt < \infty
\right\},
\end{displaymath}
equipped with the scalar product $\langle .,.\rangle_{\mathcal H}$ defined by
\begin{displaymath}
\langle\varphi,\psi\rangle_{\mathcal H} :=
\alpha_H\int_{0}^{T}\int_{0}^{T}
|t - s|^{2(H - 1)}\varphi(s)\psi(t)dsdt
\textrm{ ; }
\varphi,\psi\in\mathcal H,
\end{displaymath}
is the reproducing kernel Hilbert space of $B^H$.
%

% Proposition : Wiener integral.

%
\begin{proposition}\label{Wiener_integral}
There exists a standard Brownian motion $B$ such that :
\begin{displaymath}
B_{t}^{H} =
\int_{0}^{t}
K_H(t,s)dB_s
\textrm{ $;$ }
t\in [0,T]
\end{displaymath}
where
\begin{displaymath}
K_H(t,s) :=
c_Hs^{1/2 - H}
\int_{s}^{t}
(u - s)^{H - 3/2}u^{H - 1/2}du
\textrm{ $;$ }
(s,t)\in\Delta_T
\end{displaymath}
and $c_H > 0$ denotes a deterministic constant only depending on $H$.
\begin{displaymath}
\mathbf B^H(h) :=
\int_{0}^{T}
(K_{H}^{*}h)(t)dB_t
\textrm{ $;$ }
h\in\mathcal H
\end{displaymath}
where
\begin{displaymath}
(K_{H}^{*}h)(s) :=
\int_{s}^{T}\varphi(t)
\frac{\partial K_H}{\partial t}(t,s)dt
\textrm{ $;$ }
s\in [0,T]
\end{displaymath}
defines an iso-normal Gaussian process on $\mathcal H$ called Wiener integral with respect to $B^H$.
\end{proposition}
\noindent
That proposition summarizes several results proved at D. Nualart \cite{NUALART06}, Section 5.1.3.
\\
\\
On one hand, as an iso-normal Gaussian process, the Wiener integral defined at Proposition \ref{Wiener_integral} satisfies :
\begin{displaymath}
\forall\varphi,\psi\in\mathcal H
\textrm{, }
\mathbb E\left[\mathbf B^H(\varphi)\mathbf B^H(\psi)\right] =
\langle\varphi,\psi\rangle_{\mathcal H}.
\end{displaymath}
On the other hand, H\"older continuous functions on $[0,T]$ belong to $\mathcal H$. Then, for any (deterministic) $\beta$-H\"older continuous function $h : [0,T]\rightarrow\mathbb R$ such that $\alpha +\beta > 1$, the Young integral of $h$ with respect to $B^H$ on $[0,T]$ matches with the Wiener integral $\mathbf B^H(h)$.
\\
\\
There are many methods to simulate sample paths of a fractional Brownian motion. The most popular methods are the Wood-Chang algorithm (exact method) and the wavelet-based simulation (approximate method). Refer to T. Dieker \cite{DIEKER04} for a survey on the simulation of fractional Brownian motions.
\\
That appendix concludes on the Decreusefond-Lavaud method (cf. L. Decreusefond and N. Lavaud \cite{DL96}), particularly easy to implement. It is based on the Volterra representation of $B^H$ provided at Proposition \ref{Wiener_integral}. For $i = 0,\dots,n$, consider $t_i = iT/n$, and then
\begin{eqnarray*}
 B_{t_i}^{H} & \approx &
 \sum_{j = 0}^{i - 1}
 \left[
 \frac{1}{t_{j + 1} - t_j}
 \int_{t_j}^{t_{j + 1}}
 (t_i - t)^{H - 1/2}dt
 \right]\Delta B_{t_j}\\
 & = &
 \frac{(T/n)^{H - 1/2}}{H + 1/2}
 \sum_{j = 0}^{i - 1}
 \left[(i - j)^{H + 1/2} - (i - j - 1)^{H + 1/2}\right]\Delta B_{t_j}
\end{eqnarray*}
by putting $\Delta B_{t_j} := (T/n)^{1/2}\xi_j$ for $j = 0,\dots,n - 1$, where $\xi_0,\dots,\xi_{n - 1}$ are $n$ independent random variables of identical distributions $\mathcal N(0,1)$.
%

% Appendix B : Proofs.

%
\section{Proofs}
\noindent
At Lemma \ref{auxiliary_distributions}, the covariance function of the fractional Ornstein-Uhlenbeck process $X$ is calculated by using the construction of the reproducing kernel Hilbert space $\mathcal H$ and the Wiener integral with respect to $B^H$ defined at Appendix A, without the integration by parts formula for the Riemann integral. Proposition \ref{variance_upper_bound} allows to control, in probability, the uniform distance between the process $X$ and the solution of the associated ordinary differential equation.
\\
\\
\textit{Proof of Lemma \ref{auxiliary_distributions}.} For every $t\in\mathbb R_+$,
\begin{displaymath}
B_{t}^{H}(\vartheta) =
\mathbf B^H(\vartheta\mathbf 1_{[0,t]})
\end{displaymath}
where $\mathbf B^H$ is the Wiener integral with respect to $B^H$, defined at Proposition \ref{Wiener_integral}. Then, $B^H(\vartheta)$ is a centered Gaussian process, and for every $s,t\in\mathbb R_+$,
\begin{eqnarray*}
 R_{H,\vartheta}(s,t) & = &
 \left\langle
 \vartheta\mathbf 1_{[0,s]} ;
 \vartheta\mathbf 1_{[0,t]}\right\rangle_{\mathcal H}\\
 & = &
 \alpha_H
 \int_{0}^{s}\int_{0}^{t}
 |u - v|^{2(H - 1)}\vartheta_u\vartheta_vdudv.
\end{eqnarray*}
Since
\begin{displaymath}
X_t =
\left[C_{0}^{1-\beta} + \sigma B_{t}^{H}(\vartheta)\right]
e^{-\upsilon(1-\beta)t}
\textrm{ ; }
t\in\mathbb R_+,
\end{displaymath}
the covariance function $R_X$ satisfies :
\begin{eqnarray*}
 R_X(s,t) & = &
 \sigma^2e^{-\upsilon(1-\beta)(s + t)}
 R_{H,\vartheta}(s,t)\\
 & = &
 \alpha_H\sigma^2(1-\beta)^2
 \int_{0}^{s}\int_{0}^{t}
 |u - v|^{2(H - 1)}e^{-\upsilon(1-\beta)[(t - u) + (s - v)]}dudv
\end{eqnarray*}
for every $s,t\in\mathbb R_+$. $\square$
\\
\\
\textit{Proof of Proposition \ref{variance_upper_bound}.} For every $t\in\mathbb R_+$,
\begin{displaymath}
X_t - X_{t}^{\det} =
\upsilon(1-\beta)
\int_{0}^{t}
(X_s - X_{s}^{\det})ds +
\sigma(1-\beta)B_{t}^{H}.
\end{displaymath}
Then, $X - X^{\det}$ is an Ornstein-Uhlenbeck process, and
\begin{displaymath}
X_t - X_{t}^{\det} =
\sigma B_{t}^{H}(\vartheta)e^{-\upsilon(1-\beta)t}.
\end{displaymath}
Since $X - X^{\det}$ is a centered Gaussian process with bounded paths on $[0,T]$ ($T > 0$),
\begin{displaymath}
\forall x\in\mathbb R_{+}^{*}
\textrm{, }
\mathbb P\left(\left\|X - X^{\det}\right\|_{\infty,T} > x\right)
\leqslant
2\exp\left[-\frac{x^2}{2(\sigma^*)^2}\right]
\end{displaymath}
by Borell's inequality (cf. R.J. Adler \cite{ADLER90}, Theorem 2.1), where
\begin{eqnarray*}
 (\sigma^*)^2 & := &
 \sup_{t\in [0,T]}
 \mathbb E(|
 X_t - X_{t}^{\det}|^2)\\
 & = &
 \sigma^2R_{H,\vartheta}(T,T).
\end{eqnarray*}
That achieves the proof. $\square$
\\
\\
\textit{Proof of Corollary \ref{concentration_process_consequence}.} Consider $x > 0$ and $T > 0$. On one hand, by Proposition \ref{variance_upper_bound} :
\begin{equation}\label{inequality_consequence}
\mathbb P\left(\left\|X - X^{\det}\right\|_{\infty,T}\leqslant x\right)
\geqslant
1 - 2\exp\left[
-\frac{x^2}{2\sigma^2R_{H,\vartheta}(T,T)}\right].
\end{equation}
On the other hand, let $\omega$ be an element of $\Omega$ such that $\|X(\omega) - X^{\det}\|_{\infty,T}\leqslant x$. In other words, for every $t\in [0,T]$,
\begin{displaymath}
X_{t}^{\det} - x
\leqslant
\left|X_t(\omega)\right|
\leqslant
X_{t}^{\det} + x
\end{displaymath}
and so, by Jensen's inequality :
\begin{eqnarray*}
 0\leqslant C_t(\omega) & \leqslant & |X_{t}^{\det} + x|^{\gamma + 1}\\
 & \leqslant &
 2^{\gamma}(|X_{t}^{\det}|^{\gamma + 1} + x^{\gamma + 1})\\
 & = &
 2^{\gamma}(C_{t}^{\det} + x^{\gamma + 1}).
\end{eqnarray*}
Therefore,
\begin{displaymath}
\{\|X - X^{\det}\|_{\infty,T}\leqslant x\}
\subset
\{\forall t\in [0,T]\textrm{$,$ }
C_t\in[0,2^{\gamma}(C_{t}^{\det} + x^{\gamma + 1})]\}.
\end{displaymath}
That achieves the proof by inequality (\ref{inequality_consequence}). $\square$
\\
\\
\textit{Proof of Proposition \ref{ergodic_theorem_OU}.} The ergodic theorem provided by A. Neuenkirch and S. Tindel at \cite{NT11}, Proposition 2.3 allows to conclude if $f$ is in addition continuously differentiable. However, in the particular case of the fractional Ornstein-Uhlenbeck process, let show that the condition of Proposition \ref{ergodic_theorem_OU} is sufficient.
\\
\\
Since $Y$ is a centered, stationary and ergodic Gaussian process (cf. P. Cheridito et al. \cite{CKM03}), by the Birkhoff-Chintchin ergodic theorem together with the Fernique theorem :
\begin{displaymath}
\frac{1}{T}
\int_{0}^{T}f(Y_t)dt
\xrightarrow[T\rightarrow\infty]{\textrm{a.s.}}
\mathbb E\left[f(Y_0)\right] < \infty.
\end{displaymath}
In order to conclude, it is sufficient to show that
\begin{eqnarray*}
 \delta_T & := &
 \frac{1}{T}\left|
 \int_{0}^{T}
 f(X_t)dt -
 \int_{0}^{T}
 f(Y_t)dt\right|\\
 & &
 \xrightarrow[T\rightarrow\infty]{\textrm{a.s.}} 0,
\end{eqnarray*}
because $X_t = Y_t + (X_0 - Y_0)e^{-\upsilon(1-\beta)t}$ ; $t\in\mathbb R_+$.
\\
\\
For $T > 0$ arbitrarily chosen :
\begin{eqnarray}
 \delta_T & \leqslant &
 \frac{1}{T}
 \int_{0}^{T}
 \left|f\left[Y_t + (X_0 - Y_0)e^{-\upsilon(1-\beta)t}\right] - f(Y_t)\right|dt
 \nonumber\\
 & \leqslant &
 \label{ergodic_OU_1}
 \frac{1}{T}
 \sum_{i = 1}^{n}
 |X_0 - Y_0|^{a_i}\int_{0}^{T}
 p_i(Y_t)e^{-a_i\upsilon(1-\beta)t}dt
\end{eqnarray}
where $p_i(x) := c_i(1 + |x|)^{b_i}$ ; $x\in\mathbb R$, $i\in\{1,\dots,n\}$.
\\
\\
For $i = 1,\dots,n$, by the Cauchy-Schwarz inequality :
\begin{eqnarray*}
 \frac{1}{T}\int_{0}^{T}
 p_i(Y_t)e^{-a_i\upsilon(1-\beta)t}dt
 & \leqslant &
 \left[\frac{1}{T}\int_{0}^{T}p_{i}^{2}(Y_t)dt\right]^{1/2}
 \left[\frac{1}{T}\int_{0}^{T}e^{-2a_i\upsilon(1-\beta)t}dt\right]^{1/2}\\
 & &
 \xrightarrow[T\rightarrow\infty]{\textrm{a.s.}} 0,
\end{eqnarray*}
because
\begin{eqnarray*}
 \frac{1}{T}\int_{0}^{T}
 e^{-2a_i\upsilon(1-\beta)t}dt & = &
 -\frac{1}{2a_i\upsilon(1-\beta)T}\left[e^{-2a_i\upsilon(1-\beta)T} - 1\right]\\
 & &
 \xrightarrow[T\rightarrow\infty]{} 0
\end{eqnarray*}
and
\begin{displaymath}
\frac{1}{T}\int_{0}^{T}p_{i}^{2}(Y_t)dt
\xrightarrow[T\rightarrow\infty]{\textrm{a.s.}}
\mathbb E\left[p_{i}^{2}(Y_0)\right] < \infty
\end{displaymath}
by the Birkhoff-Chintchin ergodic theorem together with the Fernique theorem. $\square$
\\
\\
\textit{Proof of Proposition \ref{weak_upsilon_estimator}.} At the first step, it is shown that
\begin{displaymath}
\lim_{n\rightarrow\infty}
\mathbb E\left[\left|
\frac{1}{n\delta_n}
\int_{0}^{n\delta_n}
X_{t}^{2}dt -
\frac{1}{n}
\sum_{k = 0}^{n - 1}
X_{k\delta_n}^{2}
\right|\right] = 0.
\end{displaymath}
By using propositions \ref{ergodic_theorem_OU}, \ref{H_estimator} and \ref{sigma_estimator} together with the first step, the weak consistency of $\widehat\upsilon_{n}^{*}$ is shown at the second step.
\\
\\
Without loss of generality, assume that $\delta_n\in ]0,1[$ in the sequel.
\\
\\
\textbf{Step 1.} On one hand, the fractional Ornstein-Uhlenbeck process $X$ satisfies :
\begin{displaymath}
\forall p > 0\textrm{, }
\sup_{T\in\mathbb R_+}\mathbb E(|X_T|^p) <\infty.
\end{displaymath}
On the other hand, for every $k\in\{0,\dots,n - 1\}$ and $t\in [t_k,t_{k + 1}]$,
\begin{eqnarray*}
 \mathbb E^{1/2}(|X_t - X_{t_k}|^2)
 & \leqslant &
 \upsilon(1-\beta)
 \sup_{T\in\mathbb R_+}\mathbb E^{1/2}(X_{T}^{2})|t - t_k| +
 \sigma(1-\beta)
 \mathbb E^{1/2}(|B_{t}^{H} - B_{t_k}^{H}|^2)\\
 & \leqslant &
 (1-\beta)\left[
 \upsilon
 \sup_{T\in\mathbb R_+}\mathbb E^{1/2}(X_{T}^{2}) +
 \sigma\right]|t - t_k|^H.
\end{eqnarray*}
Therefore,
\begin{eqnarray*}
 \mathbb E\left[\left|
 \frac{1}{n\delta_n}
 \int_{0}^{n\delta_n}
 X_{t}^{2}dt -
 \frac{1}{n}
 \sum_{k = 0}^{n - 1}
 X_{k\delta_n}^{2}
 \right|\right] & = &
 \frac{1}{n\delta_n}\mathbb E\left[\left|
 \sum_{k = 0}^{n - 1}
 \int_{t_k}^{t_{k + 1}}
 (X_{t}^{2} - X_{t_k}^{2})dt\right|\right]\\
 & \leqslant &
 \frac{2}{n\delta_n}
 \sup_{T\in\mathbb R_+}
 \mathbb E^{1/2}(X_{T}^{2})\times\\
 & &
 \sum_{k = 0}^{n - 1}
 \int_{t_k}^{t_{k + 1}}
 \mathbb E^{1/2}(|X_t - X_{t_k}|^2)dt\\
 & \leqslant &
 \frac{C}{n\delta_n}
 \sum_{k = 0}^{n - 1}
 \int_{t_k}^{t_{k + 1}}
 |t - t_k|^Hdt\\
 & \leqslant &
 \frac{C}{H + 1}\delta_{n}^{H}
\end{eqnarray*}
with
\begin{displaymath}
C :=
2(1-\beta)\sup_{T\in\mathbb R_+}
\mathbb E^{1/2}(X_{T}^{2})
\left[
\upsilon
\sup_{T\in\mathbb R_+}\mathbb E^{1/2}(X_{T}^{2}) +
\sigma\right].
\end{displaymath}
That achieves the first step of the proof because $\delta_n\rightarrow 0$ when $n$ goes to infinity.
\\
\\
\textbf{Step 2.} Let $f :\mathbb R_{+}^{*}\times ]1/2,1[\times\mathbb R_+\rightarrow\mathbb R$ be the continuous map defined by
\begin{displaymath}
f(u,v,w) :=
\frac{1}{1-\beta}\left[\frac{w}{v^2(1-\beta)^2u\Gamma(2u)}\right]^{-1/(2u)}
\end{displaymath}
for every $u\in\mathbb R_{+}^{*}$, $v\in ]1/2,1[$ and $w\in\mathbb R_+$. By propositions \ref{H_estimator} and \ref{sigma_estimator}, and since
\begin{displaymath}
\frac{1}{n}\sum_{k = 0}^{n - 1}X_{k\delta_n}^{2}
\xrightarrow[n\rightarrow\infty]{\mathbb P}
\sigma^2(1-\beta)^{2 - 2H}\upsilon^{-2H}H\Gamma(2H)
\end{displaymath}
by the first step :
\begin{displaymath}
\widehat{\upsilon}_{n}^{*} =
f\left(\widehat{H}_n,\widehat{\sigma}_n,
\frac{1}{n}\sum_{k = 0}^{n - 1}X_{k\delta_n}^{2}\right)
\xrightarrow[n\rightarrow\infty]{\mathbb P}
\upsilon.
\end{displaymath}
That achieves the proof. $\square$
%

% References.

%

%
\end{document}